\documentclass[a4paper,12pt]{article}
\usepackage{amsmath,amscd,amssymb,french}
\newcommand{\qed}{\hfill$\blacksquare$}
\begin{document}
$\ $
\vspace{0.5cm}\\
\centerline{{\Large{\textbf{Structures de Contact sur les Vari\'et\'es 
Toriques}}}}
\vspace{0.5cm}\\
\centerline{St\'ephane DRUEL}
\begin{center}
DMI-\'Ecole Normale Sup\'erieure\\
45 rue d'Ulm\\
75005 PARIS\\
e-mail: \texttt{druel@clipper.ens.fr}
\end{center}
\vspace{1cm}
\centerline{\textbf{\S1 Introduction}}
\vspace{0.5cm}
\indent Une \textit{structure de contact} sur une vari\'et\'e
alg\'ebrique lisse est la donn\'ee d'un sous-fibr\'e
$D\subset\mathcal{T}_{X}$ de rang dim$X-1$ de sorte que la forme
$\mathcal{O}_{X}$-bilin\'eaire sur $D$ \`a valeurs dans le fibr\'e en droites
$L=\mathcal{T}_{X}/D$ d\'eduite du crochet de Lie sur
$\mathcal{T}_{X}$ soit non d\'eg\'en\'er\'ee en tout point de
$X$. Cela entra\^{\i}ne que $X$ est de dimension impaire $2n+1$ et que le
fibr\'e canonique $K_{X}$ est isomorphe \`a
$L^{-1-n}$. On peut aussi d\'efinir la structure de contact par la
donn\'ee d'un \'el\'ement $\theta\in
H^{0}(X,\Omega_{X}^{1}\otimes L)$, la \textit{forme de contact},
tel que $\theta\wedge(d\theta)^{n}$ soit partout non nul.\\
\indent Les espaces projectifs complexes et les vari\'et\'es
$\mathbb{P}_{Y}(\mathcal{T}_{Y})$, o\`u $Y$ est une vari\'et\'e lisse,
sont des exemples de vari\'et\'es de
contact.\\
\indent En quelques mots, une \textit{vari\'et\'e torique} est une
compactification \'equivariante d'un tore ([D], [O], [F]).\\
\indent L'\'etude des vari\'et\'es de contact se r\'ev\`ele
difficile. En dimension 3, elles ont \'et\'e classif\'ees par Y-G.Ye
([Y]) en utilisant la th\'eorie de Mori. D'autres auteurs \'etudient
les structures de contact sur les vari\'et\'es de Fano ([B], [L]), mais
les r\'esultats ne sont ici encore que partiels.\\
\indent Le r\'esultat principal de ce travail est le\\
\newline 
\textbf{Th\'eor\`eme}\textit{ Soit $X$ une vari\'et\'e torique
projective lisse de dimension $2n+1$ ($n\ge 1$), d\'efinie sur le corps
$\mathbb{C}$ des nombres complexes et munie
d'une structure de contact.\\
\indent Alors $X$ est soit isomorphe \`a l'espace
projectif complexe $\mathbb{P}^{2n+1}$, soit isomorphe \`a la vari\'et\'e
$\mathbb{P}_{\mathbb{P}^{1}\times\cdots\times\mathbb{P}^{1}}(\mathcal{T}_{\mathbb{P}^{1}\times\cdots\times\mathbb{P}^{1}})$.}\\
\newline
\textbf{Remerciements} : je tiens \`a exprimer toute ma gratitude \`a
A.Beauville pour m'avoir soumis ce probl\`eme et pour l'aide qu'il m'a
apport\'e.
\vspace{0.5cm}\\
\centerline{\textbf{\S2 Rappels et notations}}
\vspace{0.5cm}\\
\indent Soit $X$ une vari\'et\'e projective lisse sur le corps
$\mathbb{C}$ des nombres complexes. Le produit d'intersection entre
1-cycles et diviseurs met en dualit\'e les deux espaces vectoriels
r\'eels :
$$N_{1}(X)=(\{\text{1-cycles}\}/\equiv)\otimes\mathbb{R}\text{ et } 
N^{1}(X)=(\{\text{diviseurs}\}/\equiv)\otimes\mathbb{R},$$
\noindent o\`u $\equiv$ d\'esigne l'\'equivalence num\'erique. La
dimension commune de ces espaces vectoriels est appel\'ee le 
\textit{nombre de
Picard} de $X$. On consid\`ere le c\^one $NE(X)\subset N_{1}(X)$
engendr\'e par les classes des 1-cycles effectifs. Une \textit{raie
extr\'emale} est une demi-droite $R$ dans $\overline{NE}(X)$, 
adh\'erence de $NE(X)$ dans $N_{1}(X)$, v\'erifiant $K_{X}.R^{*}<0$  et telle que pour tout $Z_{1},Z_{2}\in\overline{NE}(X)$, si
$Z_{1}+Z_{2}\in R$ alors $Z_{1},Z_{2}\in R$ ([M]).
 Une \textit{courbe rationnelle extr\'emale} est une 
courbe rationnelle irr\'eductible $C$ telle que
$\mathbb{R}^{+}[C]$ soit une raie extr\'emale et
$-K_{X}.C\le\text{dim}X+1$. Le premier r\'esultat de la th\'eorie de Mori est que \textit{toute
raie extr\'emale est engendr\'ee par une courbe rationnelle
extr\'emale.} Le second r\'esultat fondamental est que 
\textit{toute raie 
ext\'emale $R$ admet une contraction}, c'est-\`a-dire qu'il existe une
vari\'et\'e projective normale $Y$ et un morphisme
$X\overset{\phi}{\longrightarrow}Y$, surjectif \`a fibres connexes, 
contractant les courbes irr\'eductibles $C$ telles que $[C]\in R$
(th\'eor\`eme de Kawamata-Shokurov).\\ 
\indent Dans le cas des vari\'et\'es toriques, nous avons un
r\'esultat plus pr\'ecis ([R]). Soit $X$ une vari\'et\'e torique projective
lisse associ\'ee \`a l'\'eventail $\Delta$. Nous utilisons les
notations de T.Oda ([O2]). Il existe alors des
\'el\'ements $\tau_{1},\ldots,\tau_{s}\in\Delta(d-1)$ tels que :
$$NE(X)=\overline{NE}(X)=\mathbb{R}^{+}[V(\tau_{1})]+\ldots+\mathbb{R}^{+}[V(\tau_{s})].$$
Les demi-droites $\mathbb{R}^{+}[V(\tau_{i})]$ sont appel\'ees
\textit{raies
extr\'emales g\'en\'eralis\'ees.}
Soit $R\subset NE(X)$ une raie extr\'emale g\'en\'eralis\'ee. Alors il
existe une vari\'et\'e torique projective $Y$ et un morphisme
\'equivariant $X\overset{\phi}{\longrightarrow}Y$ qui soit une
contraction au sens de Mori. Notons $A\subset X$ le lieu
exceptionnel de $\phi$ et $B=\phi(X)$ :
\begin{equation*}
\begin{CD}
X @)\phi)) Y \\
\cup & & \cup \\
A @){\phi}_{|A})) B
\end{CD}
\end{equation*}
Alors $A$ et $B$ sont deux strates toriques irr\'eductibles et le
morphisme $A\overset{\phi_{|A}}{\longrightarrow}B$ est plat et ses
fibres sont des espaces projectifs pond\'er\'es. Enfin, lorsque la
contraction est de type fibr\'ee, la vari\'et\'e $Y$ est r\'eguli\`ere et le
morphisme $\phi$ est lisse.\\
\indent Rappelons enfin un r\'esultat g\'en\'eral de 
J.Wisniewski ([W1], [W2]) sur le 
lieu exceptionnel d'une contraction extr\'emale. Soit $F$ une
composante irr\'eductible d'une fibre non triviale d'une contraction
\'el\'ementaire associ\'ee \`a la raie extr\'emale $R$. Nous appelons
\textit{lieu de $R$}, le lieu des courbes dont la classe d'\'equivalence
num\'erique appartient \`a $R$. On a alors l'in\'egalit\'e :
$$\text{dim}F+\text{dim}(\text{lieu de }R)\ge\text{dim}X+\ell(R)-1,$$
\noindent o\`u $\ell(R)$ d\'esigne la \textit{longueur} de la raie extr\'emale
$R$ :
$$\ell(R)=\text{inf}\{-K_{X}.C_{0}| C_{0} \text{ \'etant une courbe
rationnelle et }C_{0}\in R\}.$$
\vspace{0.2cm}\\
\centerline{\textbf{\S3 Preuve du th\'eor\`eme}}
\vspace{0.5cm}\\
\textbf{Lemme 1}\textit{ Soit $Y$ une vari\'et\'e torique
projective lisse de dimension $n$ et $\mathcal{E}$ un fibr\'e vectoriel
de rang $r+1$ sur $Y$. On suppose que $X=\mathbb{P}_{Y}(\mathcal{E})$
est une vari\'et\'e torique et
que le morphisme naturel
$X\overset{\phi}{\longrightarrow}Y$ est \'equivariant.\\
\indent Alors le fibr\'e $\mathcal{E}$ est totalement d\'ecompos\'e.}\\
\newline
\textit{D\'emonstration }Notons $T$ le tore de dimension $n+r$
agissant sur $X$, $T^{'}$ le tore de dimension $n$ agissant sur $Y$ et
$T\overset{\phi_{*}}{\longrightarrow}T^{'}$ le morphisme de groupes
alg\'ebriques associ\'e \`a $\phi$. Les hypoth\`eses faites entra\^{\i}nent
que ce morphisme est surjectif. Notons $N$ le noyau
$N=\ker(\phi_{*})$, c'est un groupe
alg\'ebrique produit d'un tore $T^{''}$ de dimension $r$ par un groupe
fini.\\
\indent Puisque $X$ est une vari\'et\'e torique, le fibr\'e
$\mathcal{O}_{X}(1)$ est $T$-lin\'earis\'e et en particulier
$T^{''}$-lin\'earis\'e. Le tore $T^{''}$ agissant trivialement sur
$Y$, il en r\'esulte facilement que le fibr\'e $\mathcal{E}$ est
$T^{''}$-lin\'earis\'e. Par suite, le fibr\'e $\mathcal{E}$ est somme
directe de ses composantes isotypiques :
$$\mathcal{E}=\oplus_{\chi\in X(T^{''})}\mathcal{E}_{\chi}$$
o\`u $X(T^{''})$ d\'esigne le groupe des caract\`eres de $T^{''}$. 
Si l'un des fibr\'es vectoriels $\mathcal{E}_{\chi}$ est de rang
sup\'erieur o\`u \'egal \`a deux, nous allons voir qu'un sous-tore non
trivial de
$T^{''}$ agit trivialement sur $X$, ce qui est impossible. Notons
$\chi_{1},\ldots,\chi_{k}$ $(k\ge 1)$ les caract\`eres de $T^{''}$
intervenant dans la d\'ecomposition pr\'ec\'edente et supposons 
que $k\le r-1$. Consid\'erons la composante neutre de l'intersection
des noyaux de ces caract\`eres. Par construction, l'action de ce tore
sur le fibr\'e $\mathcal{E}$ est triviale, ce qui est la conclusion
souhait\'ee.\qed\\
\newline
\textbf{Lemme 2}\textit{ Soit $X$ une vari\'et\'e torique projective lisse
de dimension $n\ge 1$ dont le fibr\'e tangent est totalement d\'ecompos\'e.\\
\indent Alors $X$ est isomorphe \`a
$\mathbb{P}^{1}\times\ldots\times\mathbb{P}^{1}$.}\\
\newline
\textit{D\'emonstration} Notons $\Delta$ l'\'eventail d\'efinissant la
vari\'et\'e torique $X$, $T$ le tore agissant sur $X$ et
$N=\text{Hom}_{gr.\,alg.}(\mathbb{G}_{m},T)$. Par hypoth\`ese, il
existe des fibr\'es en droites $L_{1},\ldots,L_{n}$ sur $X$ tels que
$\mathcal{T}_{X}=L_{1}\oplus\cdots\oplus L_{n}$. Remarquons que le fibr\'e
tangent est naturellement lin\'earis\'e puisque, si $t\in T$, $d(t.)$
d\'efinit un automorphisme du fibr\'e tangent \`a $X$. Par un r\'esultat de
A.A.Klyachko ([K] thm.1.2.3), on peut toujours supposer que l'isomorphisme
pr\'ec\'edent est \'equivariant, pour un choix convenable de
lin\'earisations des fibr\'es $L_{i}\,,\,(1\le i\le n)$.\\
\indent Prenons $\sigma\in\Delta(n)$ et consid\'erons l'ouvert
$U_{\sigma}$ qui lui est associ\'e. Puisque la
vari\'et\'e torique $X$ est lisse, $\sigma=\mathbb{R}^{+}e_{1,\sigma}+\ldots
+\mathbb{R}^{+}e_{n,\sigma}$ o\`u $(e_{1,\sigma},\ldots,e_{n,\sigma})$ est une base de $N$.
Par d\'efinition
$U_{\sigma}=\text{Spec}(\mathbb{C}[M\cap\sigma^{*}])$, o\`u $M=\text{Hom}_{gr.\,alg.}(T,\mathbb{G}_{m})=\text{Hom}_{\mathbb{Z}}(N,\mathbb{Z})$. Notons
$(x_{1,\sigma},\cdots,x_{n,\sigma})$ la base duale de
$(e_{1,\sigma},\ldots,e_{n,\sigma})$. Les $(x_{i,\sigma})_{1\le i \le n}$
sont naturellement des caract\`eres de $T$ et forment un syst\`eme de
coordonn\'ees locales sur $U_{\sigma}$. Nous notons
$\varphi_{\sigma}$ l'isomorphisme sur $\mathbb{C}^{n}$ ainsi
obtenu. L'action du tore est alors donn\'e par la formule
$t.(x_{1,\sigma},\ldots,x_{n,\sigma})=(x_{1,\sigma}(t)x_{1,\sigma},\cdots,x_{n,\sigma}(t)x_{n,\sigma})$. Aussi,
les sections
$s_{i,\sigma}=\varphi_{\sigma}^{*}(\partial_{x_{i,\sigma}})\in
H^{0}(U_{\sigma},\mathcal{T}_{X})$ sont semi-invariantes,
engendrent le fibr\'e tangent en tout point de $U_{\sigma}$ et sont
respectivement associ\'ees aux caract\`eres $(x_{i,\sigma})_{1\le i\le
n}$. On v\'erifie enfin qu'une section du fibr\'e tangent
semi-invariante relativement \`a l'un des caract\`eres $x_{i,\sigma}$ est
multiple de la section $s_{i,\sigma}$ correspondante. On peut
trouver pour chaque fibr\'e $L_{i}$, une section qui soit
semi-invariante sout $T$ et qui engendre le fibr\'e en tout point de
$U_{\sigma}$. La somme directe de ces sections fournit donc une
trivialisation du fibr\'e $\mathcal{T}_{X}$ au dessus de l'ouvert
$U_{\sigma}$ qui est d\'efinie par $n$ sections semi-invariantes.
En \'evaluant ces sections en l'unique point fixe de
$U_{\sigma}$, on v\'erifie que les caract\`eres correspondants sont
les $(x_{i,\sigma})_{1\le i\le n}$. Il en r\'esulte que
chaque section $s_{i,\sigma}$ trivialise l'un des $L_{j}$. 
Par suite, l'image de la droite
$\mathbb{C}^{*}{s_{i,\sigma}}_{|U_{\sigma\sigma'}}$ $(1\le i\le n)$ par
la fonction de transition $\varphi_{\sigma'}\varphi_{\sigma}^{-1}$ est
l'une des droites $\mathbb{C}^{*}{s_{i,\sigma'}}_{|U_{\sigma\sigma'}}$
$(1\le i\le n)$. Autrement dit, la matrice associ\'ee \`a la fonction de
transition $\varphi_{\sigma'}\varphi_{\sigma}^{-1}$ exprim\'ee dans les
bases $({s_{i,\sigma'}}_{|U_{\sigma\sigma'}})_{1\le i\le n}$ et
$({s_{i,\sigma}}_{|U_{\sigma\sigma'}})_{1\le i\le n}$ n'a qu'un seul
terme non nul par ligne et par colonne.\\
\indent C'est ce r\'esultat que nous allons maintenant exploiter.
Consid\'erons un c\^one $\tau\in\Delta(n-1)$ et soit $\sigma$
et $\sigma'$ les deux c\^ones de dimension $n$ s\'epar\'es par
$\tau$. Il existe donc une famille $(e_{1},\ldots,e_{n-1})$ d'\'el\'ements
de $N$ et deux autres \'el\'ements de $N$, $e_{n}$ et $e_{n+1}$, tels
que les familles $(e_{1},\cdots,e_{n-1},e_{n})$ et
$(e_{1},\ldots,e_{n-1},e_{n+1})$ soient des bases de $N$ et tels que
l'on ait $\tau=\mathbb{R}^{+}e_{1}+\cdots+\mathbb{R}^{+}e_{n-1}$,
$\sigma=\tau+\mathbb{R}^{+}e_{n}$ et
$\sigma'=\tau+\mathbb{R}^{+}e_{n+1}$. On a en outre une relation :
$$e_{n+1}+e_{n}+\sum_{i=1}^{n-1}\alpha_{i}e_{i}=0,$$
\noindent o\`u $\alpha_{i}\in\mathbb{Z}$.
Notons $(x_{1},\cdots,x_{n})$ la base duale de la base
$(e_{1},\ldots,e_{n-1},e_{n})$ et $(z_{1},\cdots,z_{n})$ la base duale
de la base $(e_{1},\ldots,e_{n-1},e_{n+1})$. Le changement de cartes
est alors donn\'e par la formule :
$$(z_{1},\cdots,z_{n})\mapsto
(\frac{z_{1}}{z_{n}^{\alpha_{1}}},\cdots,\frac{z_{n-1}}{z_{n}^{\alpha_{n-1}}},\frac{1}{z_{n}}).$$
La fonction de transition $\varphi_{\sigma\sigma'}$ est donc donn\'ee par la matrice :
$$
\left(
\begin{array}{cccccc}
z_{n}^{-\alpha_{1}} & 0 & \ldots & \ldots & 0 &
-\alpha_{1}z_{1}z_{n}^{-\alpha_{1}-1} \\
0 & z_{n}^{-\alpha_{2}} & \ddots & & \vdots & \vdots \\
\vdots & \ddots & \ddots & \ddots & \vdots & \vdots \\
\vdots & & \ddots & \ddots & 0 & \vdots \\
\vdots & & & \ddots & z_{n}^{-\alpha_{n-1}} &
-\alpha_{n-1}z_{1}z_{n}^{-\alpha_{n-1}-1} \\
0 & \ldots & \ldots & \ldots & 0 & -z_{n}^{-2} \\
\end{array}
\right)
$$
Il en r\'esulte les \'egalit\'es
$\alpha_{1}=\ldots=\alpha_{n-1}=0$ et donc $e_{n+1}=-e_{n}$.\\
\indent Soit $e\in\Delta(1)$ et soit $\sigma\in\Delta(n)$ un c\^one
dont $e$ est une face de dimension 1. On identifie le c\^one $e$ et
l'\'el\'ement primitif de $N$ qui le d\'etermine. Alors il existe
une base $(e_{1}=e,\ldots,e_{n})$ de $N$ telle que
$\sigma=\mathbb{R}^{+}e_{1}+\ldots+\mathbb{R}^{+}e_{n}$. Le c\^one
$\mathbb{R}^{+}e_{2}+\ldots+\mathbb{R}^{+}e_{n}$ est une face de
dimension $n-1$ de $\sigma$ et donc d'apr\`es ce qui pr\'ec\`ede, on a
que $-e\in \Delta(1)$ et que le c\^one
$\mathbb{R}^{+}(-e_{1})+\ldots+\mathbb{R}^{+}e_{n}\in\Delta(n-1)$. Il
en r\'esulte que $\pm e_{i}\in\Delta(1)$ $(1\le i\le n)$ et que tous
les c\^ones possibles de dimension $n$ form\'es sur ces vecteurs sont
des c\^ones de l'\'eventail $\Delta$. Si $e'\in\Delta(1)$, alors $e'$
est contenu dans l'un des c\^ones pr\'ec\'edents et donc $e'$ est l'un
des $(\pm e_{i})_{1\le i\le n}$ et on a donc d\'etermin\'e tous les
c\^ones de dimension maximale,
ce qui d\'etermine $\Delta$ et termine la preuve du lemme.\qed\\
\newline
\textbf{Th\'eor\`eme}\textit{ Soit $X$ une vari\'et\'e torique
projective lisse de dimension $2n+1$ ($n\ge 1$) munie d'une structure
de contact.\\   
\indent Alors $X$ est soit isomorphe \`a l'espace
projectif complexe $\mathbb{P}^{2n+1}$, soit isomorphe \`a la
vari\'et\'e 
$\mathbb{P}_{\mathbb{P}^{1}\times\cdots\times\mathbb{P}^{1}}(\mathcal{T}_{\mathbb{P}^{1}\times\cdots\times\mathbb{P}^{1}})$.}\\
\newline
\textit{D\'emonstration} La preuve de ce th\'eor\`eme repose sur
l'\'etude des contractions extr\'emales de $X$. Notons
$T$ le tore agissant sur $X$. Puisque le
diviseur canonique $K_{X}$ de $X$ n'est pas effectif, il existe une
raie extr\'emale $R$ au sens de Mori engendr\'ee par une courbe
rationnelle lisse $C$ $T$-invariante. Puisque $K_{X}=-(n+1)L$, on a
$\ell(R)=n+1$ ou $\ell(R)=2n+2$, o\`u $\ell(R)$ d\'esigne la longueur
de la raie extr\'emale $R$. Dans ce dernier cas, $X$ est de Fano et
$\text{Pic}(X)=\mathbb{Z}$ ([W1] prop. 2.4). Il en r\'esulte que $X$
est isomorphe \`a l'espace projectif complexe $\mathbb{P}^{2n+1}$
([O1] thm.7.1).\\
\indent Il nous reste donc \`a \'etudier le cas o\`u $\ell(R)=n+1$. Notons
$X\overset{\phi}{\longrightarrow}Y$ la contraction de Mori associ\'ee
\`a la raie $R$, o\`u $Y$ est une vari\'et\'e torique projective. 
Notons $T^{'}$ le tore agissant sur $Y$. Notons $A\subset X$ le lieu
exceptionnel de $\phi$ et $B=\phi(X)$ :
\begin{equation*}
\begin{CD}
X @)\phi)) Y \\
\cup & & \cup \\
A @){\phi}_{|A})) B
\end{CD}
\end{equation*}
Puisque la courbe $C$ est $T$-invariante, elle est contract\'ee sur un
point fixe $y\in B$ sous $T^{'}$. Le morphisme
$\phi$ \'etant \'equivariant, la fibre $F=\phi^{-1}(y)$ au dessus
de $y$ est donc $T$-invariante. Par suite,
puisque cette fibre est irr\'eductible, on en d\'eduit que $F$ est une
strate torique et donc que $F$ est lisse puisque $X$ l'est. La fibre
$F$ est donc isomorphe \`a un espace
projectif $\mathbb{P}^{k}$ ($k\ge 1)$ et, puisque $C$ est une strate
torique de la vari\'et\'e torique $\mathbb{P}^{k}$, il en r\'esulte
que la courbe
$C$ s'identifie \`a une droite de $\mathbb{P}^{k}$. Nous
savons aussi qu'il existe une courbe rationnelle $C_{1}\subset X$
telle que $C_{1}\in R$ et $-K_{X}.C_{1}=n+1$. Mais, puisqu'il existe
un diviseur $D$ sur $X$ tel que $D.C=1$ et puisque
$R=\mathbb{R}^{+}[C]$, il en r\'esulte facilement que $C\equiv
C_{1}$, de sorte que $C.L=1$. Par suite, on a
$(F,L_{|F})=(\mathbb{P}^{k},\mathcal{O}_{\mathbb{P}^{k}}(1))$.\\
\indent Puisque
$H^{0}(\mathbb{P}^{k},\Omega_{\mathbb{P}^{k}}^{1}(1))=(0)$, la
restriction de la forme de contact $\theta$ \`a $\mathbb{P}^{k}$ est
nulle et on v\'erifie par un calcul
en coordonn\'ees locales que, pour tout $x\in\mathbb{P}^{k}$, le sous
espace vectoriel
$\mathcal{T}_{\mathbb{P}^{k}}(x)\subset D(x)$ est totalement
isotrope pour la forme altern\'ee de contact qui, par hypoth\`ese est
non-d\'eg\'en\'er\'ee. Il en r\'esulte l'in\'egalit\'e $k\le n$  
et donc l'\'egalit\'e $k=n$ puisque
$\text{dim}F\ge\ell(R)-1=n$. Utilisant \`a nouveau l'in\'egalit\'e
 de Wisniewski, on en d\'eduit que la
contraction est de type fibr\'ee et donc que $Y$ est r\'eguli\`ere, de
dimension $n+1$ et que $\phi$ est lisse.\\
\indent  Consid\'erons la
suite exacte :
$$0\longrightarrow\mathcal{T}_{X/Y}\longrightarrow\mathcal{T}_{X} 
\longrightarrow\phi^{*}\mathcal{T}_{Y}\longrightarrow 0$$
La fl\`eche $\mathcal{T}_{X/Y}\longrightarrow L$ obtenue par
composition avec la projection $\mathcal{T}_{X}\longrightarrow L$
\'etant identiquement nulle par le th\'eor\`eme de Grauert et les
r\'esultats ci-dessus, il existe une fl\`eche surjective
$\phi^{*}\mathcal{T}_{Y}\longrightarrow L \longrightarrow 0$ et donc
un morphisme $X\longrightarrow \mathbb{P}_{Y}(\mathcal{T}_{Y})$ au
dessus de $Y$ qui induit un isomorphisme sur chaque fibre. Il en
r\'esulte que ce morphisme est en fait un
isomorphisme. Le r\'esultat d\'ecoule alors des lemmes 1 et 2.\qed
\newpage
\centerline{\textbf{R\'ef\'erences}}
$\ $
\newline
\noindent [B] A.Beauville, \emph{Fano Contact Manifolds and Nilpotent
Orbits}, Comment. Math. Helvet., \`a para\^{\i}tre.\\
\newline
[D] V.I.Danilov, \emph{The geometry of toric varieties},
Russian Math. Surveys, 33, 97-154, 1978.\\
\newline
[F] W.Fulton, \emph{Introduction to Toric Varieties}, Annals of
Mathematics Studies, Princeton
University Press, 131, 1993.\\
\newline 
[K] A.A.Klyachko, \emph{Equivariant Bundles on Toral Varieties},
Math. USSR Izv., 35, 337-375, 1990.\\
\newline
[L] C.Lebrun, \emph{Fano manifolds, contact structures and
quaternionic geometry}, Int. journ. of Math., 6, 419-437, 1995.\\
\newline
[M] S.Mori, \emph{Threefolds whose canonical bundles are not
numerically effective}, Ann. Math., 116, 133-176, 1982.\\
\newline
[O1] T.Oda, \emph{Lectures on Torus Embeddings and Applications}, Tata
Inst. of Fund. Research, 58, Springer, 1978.\\
\newline
[O2] T.Oda, \emph{Convex Bodies and Algebraic Geometry}, 15,
Springer-Verlag, 1988.\\
\newline
[R] M.Reid, \emph{Decomposition of toric morphism}, dans
\emph{Arithmetic and Geometry, papers dedicated to I.R. Shafarevitch
on the occasion of his 60th birthday}, Progress in Math. 36,
Birkha\"{u}ser, 395-418, 1983.\\
\newline
[W1] J.Wisniewski, \emph{Length of extremal rays and generalized
adjonction}, Math. Z., 200, 409-427, 1989.\\
\newline
[W2] J.Wisniewski, \emph{On contractions of extremal rays on Fano
manifolds}, J. reine u. angew Math., 417, 141-157, 1991.\\
\newline
[Y] Y-G.Ye, \emph{A note on complex projective threefolds admitting
holomorphic contact structures}, Invent. Math., 121, 421-436, 1995.

\end{document}